\documentclass[12pt]{amsart}

\usepackage{amsmath}
\usepackage{amsthm}
\usepackage{amsfonts}
\usepackage{amssymb,amsmath,latexsym}
\usepackage{graphicx}

\theoremstyle{plain}
\newtheorem{theorem}{Theorem}[section]
\newtheorem{lemma}[theorem]{Lemma}
\newtheorem{corollary}[theorem]{Corollary}

\theoremstyle{definition}
\newtheorem{remark}[theorem]{Remark}

\title[Asymptotic probability]{An application of Jacobi's elliptic functions to  asymptotic probabilities for conformal restriction measures}

\author{Robert~O. Bauer}
\address{Altgeld Hall\\Department of Mathematics\\ 
	University of Illinois at Urbana-Champaign\\ 
	1409 West Green Street \\ 
	Urbana, IL 61801, USA}
\email{rbauer@math.uiuc.edu}

\begin{document}

\begin{abstract}
We show that for the conformal restriction measure with exponent $b$ in the unit disk on hulls $\gamma$ connecting $e^{ix}$ to 1 the probability of the event that $\gamma$ avoids the disk of radius $q$ centered at zero decays like $\exp(-b\pi x/(1-q))$ if either $b\in[5/8,1]\cup[5/4,\infty)$ and $x\in(0,\pi]$, or if  $b\in(1,5/4)$, $x\in(0,\pi)$, and $bx\le\pi$. \end{abstract}

\maketitle

\section{Introduction}

For a simply connected domain $D$ and two distinct points $z_1,z_2$ on the boundary of $D$ denote $P^b_{D,z_1, z_2}$ a conformal restriction measure of exponent $b$ supported on closed, simply connected subdomains $\gamma$ of $\overline{D}$ with $\gamma\cap\partial D=\{z_1,z_2\}$, see \cite{LSW:rest}, \cite{lawler:2005}. This measure exists for all $b\ge 5/8$. For example, the  chordal Schramm-Loewner evolution with parameter $\kappa=8/3$ gives rise to the conformal restriction measure with exponent $b=5/8$, and the {\em filling} of a Brownian excursion gives rise to the conformal restriction measure with exponent $b=1$. We will assume that the boundary of $D$ near $z_1$ and $z_2$ is smooth. Let $a<0$, $q=e^a$, and denote $A_q$ the annulus $\{z:q<|z|<1\}$, and $\mathbb U$ the unit disk $\{z:|z|<1\}$. 
In this note we derive the asymptotics of the non-intersection probability
\[
	(a,b,x)\in[-\infty,0]\times[5/8,\infty)\times[0,2\pi]\mapsto F(a,b,x)\equiv P^b_{\mathbb U,e^{ix},1}(\gamma\subset A_q)
\]
as $a\nearrow0$. The defining properties for conformal restriction measures are {\em conformal invariance}, i.e. for any conformal map $f:D\to f(D)$ and any subdomain $D'$ of $D$, we have\smallskip

$\bullet\ P^b_{D,z_1, z_2}(\gamma\subset D')=P^b_{f(D),f(z_1), f(z_2)}(\gamma\subset f(D')),$\smallskip

\noindent{and} {\em conformal restriction} in the sense that if $D'$ is a simply connected subdomain of $D$ such that $D\backslash D'$ is bounded away from $z_1, z_2$, and $f$ is a conformal map from $D'$ onto $D$ such that $f(z_{1,2})=z_{1,2}$, then\smallskip

$
\bullet\ P^b_{D,z_1, z_2}(\gamma\subset D')=|f'(z_1)f'(z_2)|^{b}.
$\smallskip

\noindent{We} will also use a generalization of a result of Beffara, see \cite{beffara}, which says that if $D'$ is  a subdomain of $D$, not necessarily simply connected, such that $\partial D\subset\partial D'$, and $f$ maps $D'$ conformally onto $f(D')$ such that $f(D')\subset D$, $f(\partial D)=\partial D$, and $f(z_{1,2})=z_{1,2}$, then \smallskip

$
\bullet\ P^b_{D,z_1, z_2}(\gamma\subset D')=P^b_{D,z_1, z_2}(\gamma\subset f(D'))|f'(z_1)f'(z_2)|^{b}.
$
{\smallskip

\noindent{We}} give a proof of this fact in a particular case in Lemma \ref{L:change}.

The key to our analysis is the transformation of the annulus $A_q$ to the unit disk $\mathbb U$ with a horizontal slit $[-L,L]$ along the real axis. What our estimate shows for example is that the probability that $\text{SLE}_{8/3}$ in the unit disk from $i$ to $-i$ stays in a thin annulus $A_q$ is, to leading order, the same as hitting at least one of the two real segments $(-1,-L),(L,1)$ which are not part of the slit.

An estimate closely related to ours appears in \cite[Lemma 18]{werner:2005}. There, the general form of the estimate is derived from an excursion representation of $\text{SLE}_{8/3}$ and it is stated that, and briefly indicated how, the explicit values of the constants can be derived from a comparison argument. In this note we carry through such a comparison argument in detail, and the general form of the estimate is established together with the constants at once. The upper bound is more subtle and requires the majority of the work. The upper bound, or rather our lack of finding a better one, is also the reason why there is gap in the parameter range for which we obtain the  asymptotic behavior. 

A related estimate for the asymptotic behavior also appears in \cite{cardy:loop}, derived using Coulomb gas techniques. In both, \cite{werner:2005} and \cite{cardy:loop}, the aim is to find---asymptotically---the weight, according to the conformally invariant measure on self-avoiding loops, of the loops which surround an annulus. 

The following proof is a straightforward modification and generalization (from $b=5/8$ to $b\in[5/8,\infty)$) of the proof contained in our paper \cite{SLE83}.

\section{Asymptotic behavior of the non-intersection probability}\label{S:asymp}

For each $q\in[0,1)$ there exists a unique $L=L(q)\in[0,1)$ such that $A_q$ and $\mathbb U\backslash[-L,L]$ are conformally equivalent. As $q$ increases to 1, $L$ increases to 1 as well. Denote $f$ the conformal equivalence, normalized by $f(1)=1$. For $x\in(0,\pi]$, let $z_1=e^{ix/2}, z_2=e^{-ix/2}$. By symmetry, if $w_{1,2}=f(z_{1,2})$, then $w_2=\bar w_1$.

In what follows we will mean by $h(a)\asymp g(a)$ as $a\nearrow0$, that \[\lim_{a\nearrow0}\log h(a)/\log g(a)=1.
\]

\begin{lemma}\label{L:sn-asymptotics} For $x\in(0,\pi]$, we have
\[
1-L\asymp e^{\frac{\pi^2}{4a}},\text{ and }|f'(z_1)|\asymp|1-f(z_1)|\asymp e^{\frac{\pi}{4a}(\pi-x)}\]
as $a\nearrow0$. 
\end{lemma}

\begin{proof}
From \cite[Chap. VI, Sec. 3]{nehari:1952}, 
\[
f(z)=L\text{ sn}\left(\frac{2iK}{\pi}\log\frac{z}{q}+K;q^4\right),
\]
where $\text{sn}(z)$ is the analytic function for which $\text{sn}'(0)=1$ and which maps the rectangle $\{z:-K<\Re z<K,0<\Im z<i K'\}$ onto the upper half-plane in such a way that $\text{sn}(\pm K)=\pm1$ and $\text{sn}(\pm K+iK')=\pm k^{-1}$. Furthermore, $q^4=\exp(-\pi K'/K)$, and $L=\sqrt{k}$. It is classical that $\text{sn}'(z)=[(1-\text{sn}^2(z))(1-k^2\text{sn}^2(z))]^{1/2}$. Thus
\begin{equation}\label{E:fprime}
f'(z)=(2iK/\pi z)[(L^2-f^2(z))(1-L^2 f^2(z))]^{1/2}.
\end{equation}
Define $h, \tau$ by $q^4=h=e^{i\pi\tau}$, and set $v=\frac{i}{\pi}\log\frac{z_1}{q}+\frac{1}{2}$. Then it follows from \cite[II, 3.]{hurwitz:1964}, 
that 
\[
L=\frac{\theta_2(0|\tau)}{\theta_3(0|\tau)}, \text{ and } f(z)=\frac{\theta_1(v|\tau)}{\theta_0(v|\tau)}.
\]
Here
\begin{align}
\theta_1(v|\tau)&=2\sum_{n=0}^{\infty}(-1)^n h^{(n+1/2)^2}\sin(2n+1)\pi v,\notag\\
\theta_2(v|\tau)&=2\sum_{n=0}^{\infty}h^{(n+1/2)^2}\cos(2n+1)\pi v,\notag\\
\theta_3(v|\tau)&=1+2\sum_{n=1}^{\infty}h^{n^2}\cos 2n\pi v,\notag\\
\theta_0(v|\tau)&=1+2\sum_{n=1}^{\infty}(-1)^n h^{n^2}\cos2n\pi v.\notag
\end{align}
Using linear transformations of theta functions we may write
\[
\frac{\theta_2(0|\tau)}{\theta_3(0|\tau)}=\frac{\theta_0(0|-\frac{1}{\tau})}{\theta_3(0|-\frac{1}{\tau})},\text{ and } \frac{\theta_1(v|\tau)}{\theta_0(v|\tau)}=i\frac{\theta_1(\frac{v}{\tau}|-\frac{1}{\tau})}{\theta_2(\frac{v}{\tau}|-\frac{1}{\tau})}.
\]
Hence, if $h'=\exp(-i\pi/\tau)$, and using the series representation of $\theta_0$ and $\theta_3$, we get
\[
L=\frac{1+2\sum_{n=1}^{\infty}(-1)^n(h')^{n^2}}{1+2\sum_{n=1}^{\infty}(h')^{n^2}}=1-4h'+O((h')^2),
\]
which is the first statement of the lemma. For the second, we use the infinite product representation of $\theta_1$ and $\theta_2$, giving
\[
i\frac{\theta_1(\frac{v}{\tau}|-\frac{1}{\tau})}{\theta_2(\frac{v}{\tau}|-\frac{1}{\tau})}=\frac{e^{2i\pi v/\tau}-1}{e^{2i\pi v/\tau}+1}\prod_{n=1}^{\infty}\frac{(1-(h')^{2n}e^{2i\pi v/\tau})(1-(h')^{2n}e^{-2i\pi v/\tau})}{(1+(h')^{2n}e^{2i\pi v/\tau})(1+(h')^{2n}e^{-2i\pi v/\tau})}.
\]
Since $\exp(2i\pi v/\tau)=i\exp(-(\pi/4a)(\pi-x))$, the infinite product is $1+O(\exp(\pi^2/(4a)))$, and
\[
\frac{e^{2i\pi v/\tau}-1}{e^{2i\pi v/\tau}+1}=1+2ie^{\frac{\pi}{4a}(\pi-x)}+O(e^{\pi^2/(4a)}),
\]
as $a\nearrow0$. Using equation \eqref{E:fprime}, the lemma now follows. 
\end{proof}

Recall that $z_1=e^{ix/2}$, $w_1=f(z_1)$, and set $u=i(1+w_1)/(1-w_1)$. The following result is analogous to a result in \cite{beffara}. We will give a direct argument.

\begin{lemma}\label{L:change} The probability $P^b_{\mathbb U,e^{ix},1}(\gamma\subset A_q)$ is equal to
\[P^b_{\mathbb H,u,-u}\left(\gamma\cap i\left[\frac{1-L}{1+L},\frac{1+L}{1-L}\right]=\emptyset\right)\left|\frac{f'(z_1)(1-z_1)}{1-f(z_1)}\right|^{2b}.
\]
\end{lemma}

\begin{proof}
Denote $B$ a simple curve connecting the inner and outer boundary of $A_q$, so that $B$ is bounded away from $z_1$ and $z_2$. Denote $\phi$ a conformal map from $A_q\backslash B$ onto $\mathbb U$ such that $\phi(z_{1,2})=z_{1,2}$, and $\psi$ a conformal map from $f(A_q\backslash B)$ onto $\mathbb U$ such that $\psi(w_{1,2})=w_{1,2}$. Then, by conformal restriction, 
\begin{align}\label{E:psis}
P^b_{\mathbb U,z_1, z_2}(\gamma\subset A_q\backslash B)&=|\phi'(z_1)\phi'(z_2)|^{b},\notag\\
P^b_{\mathbb U,w_1, w_2}(\gamma\subset f(A_q\backslash B))&=|\psi'(w_2)\psi'(w_2)|^{b}. 
\end{align}
Since $T\equiv\phi\circ f\circ \psi^{-1}$ maps $\mathbb U$ onto $\mathbb U$ and sends $w_{1,2}$ to $z_{1,2}$, there is a  pair $w_0, z_0\in\partial\mathbb U$ such that $T$ is the linear transformation given by
\[
\frac{T(w)-w_1}{T(w)-w_2}\cdot\frac{w_0-w_2}{w_0-w_1}=\frac{z-z_1}{z-z_2}\cdot\frac{z_0-z_2}{z_0-z_1}.
\]
A calculation gives
\[
T'(w_1)T'(w_2)=\left(\frac{z_1-z_2}{w_1-w_2}\right)^2,
\]
which together with $|f'(z_1)|=|f'(z_2)|$ implies
\begin{equation}\label{E:part}
P^b_{\mathbb U,z_1, z_2}(\gamma\subset A_q\backslash B) =P^b_{\mathbb U,w_1, w_2}(\gamma\subset f(A_q\backslash B))\left|\frac{f'(z_1)(z_1-z_2)}{w_1-w_2}\right|^{2b}.
\end{equation}
By an inclusion/exclusion argument, equation \eqref{E:part} also holds if $A_q\backslash B$ is replaced by $A_q$. Finally, by conformal invariance, 
\[
P^b_{\mathbb U,w_1, w_2}(\gamma\subset f(A_q))=P^b_{\mathbb H,u,-u}\left(\gamma\cap i\left[\frac{1-L}{1+L},\frac{1+L}{1-L}\right]=\emptyset\right).
\]
\end{proof}

Note that because $x\in(0,\pi]$ we have $\arg z_1,\arg w_1\in(0,\pi/2]$ and so $u\le-1.$ We will use the following lower and upper bounds:
\begin{align}\label{E:lower}
&P^b_{\mathbb H,u,-u}\left(\gamma\cap i\left[\frac{1-L}{1+L},\frac{1+L}{1-L}\right]=\emptyset\right)\notag\\
&\ge P^b_{\mathbb H,u,-u}(\gamma\cap i(0,\frac{1+L}{1-L}]=\emptyset)+P^b_{\mathbb H,u,-u}(\gamma\cap i[\frac{1-L}{1+L},\infty)=\emptyset)\notag\\
&=P^b_{\mathbb H,u,-u}(\gamma\cap i(0,\frac{1+L}{1-L}]=\emptyset)+P^b_{\mathbb H,\frac{1}{u},-\frac{1}{u}}(\gamma\cap i(0,\frac{1+L}{1-L}]=\emptyset),
\end{align}
and
\begin{align}\label{E:upper}
&P^b_{\mathbb H,u,-u}\left(\gamma\cap i\left[\frac{1-L}{1+L},\frac{1+L}{1-L}\right]=\emptyset\right)\notag\\
&\le P^b_{\mathbb H,u,-u}(\gamma\cap i(0,\frac{1+L}{1-L}]=\emptyset)+P^b_{\mathbb H,\frac{1}{u},-\frac{1}{u}}(\gamma\cap i[\frac{1+L}{1-L},\infty)=\emptyset)\notag\\
&\quad+P^b_{\mathbb H,u,-u}(\gamma\cap i(0,\frac{1-L}{1+L})\neq\emptyset,\gamma\cap i(\frac{1+L}{1-L},\infty)\neq\emptyset).
\end{align}

For $c\in\mathbb R, d>0$, set 
\[
g_{c,d}(z)=\frac{|c|}{\sqrt{c^2+d^2}}\sqrt{z^2+d^2}.
\]
Then $g_{c,d}$ maps $\mathbb H\backslash i(0,d]$ conformally onto $\mathbb H$ such that $g_{c,d}(\pm c)=\pm c$. Furthermore,
\[
|g'_{c,d}(c)g'_{c,d}(-c)|=\frac{c^4}{(c^2+d^2)^2},
\]
and so by conformal restriction
\begin{equation}\label{E:gcd}
P^b_{\mathbb H,c, -c}(\gamma\cap i(0,d]=\emptyset)=[c^2/(c^2+d^2)]^{2b}.
\end{equation}

\begin{corollary}\label{C:gcd}
We have
\[
P^b_{\mathbb H,u,-u}(\gamma\cap i(0,\frac{1+L}{1-L}]=\emptyset)+P^b_{\mathbb H,\frac{1}{u},-\frac{1}{u}}(\gamma\cap i(0,\frac{1+L}{1-L}]=\emptyset)\asymp e^{\frac{b\pi x}{a}}
\]
as $a\nearrow0$. 
\end{corollary}

\begin{proof}
By \eqref{E:gcd},
\[
P^b_{\mathbb H,u,-u}(\gamma\cap i(0,\frac{1+L}{1-L}]=\emptyset)
=\left(\frac{u(1-L)}{1+L}\right)^{4b}\left(1+\frac{(u(1-L))^2}{(1+L)^2}\right)^{-2b},
\]
and from Lemma \ref{L:sn-asymptotics}
\[
\left(\frac{u(1-L)}{1+L}\right)^{4b}\left(1+\left(\frac{u(1-L)}{1+L}\right)^2\right)^{-2b}\asymp e^{\frac{b\pi x}{a}}.
\]
Similarly,
\[
P^b_{\mathbb H,\frac{1}{u},-\frac{1}{u}}(\gamma\cap i(0,\frac{1+L}{1-L}]=\emptyset)\asymp e^{\frac{b\pi^2}{a}+\frac{b\pi}{a}(\pi-x)},
\]
so that this term is negligible compared to the first if $0<x<\pi$, and of the same order if $x=\pi$.
\end{proof}

\begin{lemma}\label{L:upper-lower} We have
\[
P^b_{\mathbb H,u,-u}\left(\gamma\cap i(0,\frac{1-L}{1+L})\neq\emptyset,\gamma\cap i(\frac{1+L}{1-L},\infty)\neq\emptyset\right)\asymp e^{\pi^2/a},
\]
as $a\nearrow0$.
\end{lemma}

\begin{proof} First,
\begin{align}\label{E:first}
&P^b_{\mathbb H,u,-u}\left(\gamma\cap i(0,\frac{1-L}{1+L})\neq\emptyset,\gamma\cap i(\frac{1+L}{1-L},\infty)\neq\emptyset\right)\notag\\
&=P^b_{\mathbb H,u,-u}\left(\gamma\cap i(0,\frac{1-L}{1+L})\neq\emptyset\right)+P^b_{\mathbb H,\frac{1}{u},-\frac{1}{u}}\left(\gamma\cap i(0,\frac{1-L}{1+L})\neq\emptyset\right)\notag\\
&\quad-P^b_{\mathbb H,u,-u}\left(\gamma\cap i((0,\frac{1-L}{1+L})\cup(\frac{1+L}{1-L},\infty))\neq\emptyset\right).
\end{align}
The last probability on the right equals
\[
P^b_{\mathbb U,w_1, w_2}(\gamma\cap((-1,-L]\cup[L,1))\neq\emptyset).
\]
To calculate this probability, note that 
\[
g_L(w)\equiv\frac{1+w^2-\sqrt{(1+w^2)^2-4p^2w^2}}{2pw}
\]
maps $\mathbb U\backslash((-1,-L]\cup[L,1))$ onto $\mathbb U$ if $2p=(L+1/L)$, see \cite[Chapter 3]{ivanov}. Here, the square root is chosen so that $g_L(i)=i$. Setting $w=e^{i\varphi}$, this can be written
\begin{equation}
g_L(w)=\begin{cases}\frac{1}{p}\cos\varphi+i\sqrt{1-\frac{1}{p^2}\cos^2\varphi}, &\text{if $\varphi\in(0,\pi/2]$;}\\
\frac{1}{p}\cos\varphi-i\sqrt{1-\frac{1}{p^2}\cos^2\varphi}, &\text{if $\varphi\in[-\pi/2,0)$.}
\end{cases}
\end{equation}
Then
\[
g'_L(w)g'_L(\bar w)=-\frac{\sin^2\varphi}{p^2-1+\sin^2\varphi}.
\]
Denote $T$ a (fractional) linear transformation from $\mathbb U$ onto $\mathbb U$ sending $g_L(w_{1,2})$ onto $w_{1,2}$. Then, as in the proof of Lemma \ref{L:change}, 
\[
T'(g_L(w_1))T'(g_L(w_2))=\frac{\sin^2\varphi}{1-\frac{1}{p^2}\cos^2\varphi},
\]
where now $\varphi=\arg w_1$. Thus, by conformal restriction,
\begin{equation}\label{E:p1}
P^b_{\mathbb U,w_1, w_2}(\gamma\cap((-1,-L]\cup[L,1))\neq\emptyset)=1-\left[\frac{p\sin^2\varphi}{p^2-1+\sin^2\varphi}\right]^{2b}.
\end{equation}
Finally, from the definition of $u$ and $\varphi$ in terms of $w_1$, it follows that $u=-\cot(\varphi/2)$ and so $4/\sin^2\varphi=(u+1/u)^2$. A calculation now gives
\begin{align}\label{E:p2}
&\frac{p^2-1+\sin^2\varphi}{p\sin^2\varphi}\notag\\
&=1+\left(\frac{1-L}{1+L}\right)^2(u^2+\frac{1}{u^2})+\frac{(1-L)^4}{8(L+L^3)}\left[2+\left(\frac{1-L}{1+L}\right)^2(u^2+\frac{1}{u^2})\right].
\end{align}
On the other hand, \eqref{E:gcd} implies
\begin{equation}\label{E:hit1}
P^b_{\mathbb H,u,-u}\left(\gamma\cap i(0,\frac{1-L}{1+L})\neq\emptyset\right)=1-\left(1+\left(\frac{1-L}{1+L}\right)^2\frac{1}{u^2}\right)^{-2b}
\end{equation}
and
\begin{equation}\label{E:hit2}
P^b_{\mathbb H,\frac{1}{u},-\frac{1}{u}}\left(\gamma\cap i(0,\frac{1-L}{1+L})\neq\emptyset\right)=1-\left(1+\left(\frac{1-L}{1+L}\right)^2 u^2\right)^{-2b}.
\end{equation}
Combining \eqref{E:hit1}, \eqref{E:hit2}, \eqref{E:p1}, and \eqref{E:first}, we get
\begin{align}\label{E:ende}
&P^b_{\mathbb H,u,-u}\left(\gamma\cap i(0,\frac{1-L}{1+L})\neq\emptyset,\gamma\cap i(\frac{1+L}{1-L},\infty)\neq\emptyset\right)\notag\\
&=1-\left(1+\left(\frac{1-L}{1+L}\right)^2\frac{1}{u^2}\right)^{-2b}
+1-\left(1+\left(\frac{1-L}{1+L}\right)^2 u^2\right)^{-2b}\notag\\
&\quad-1+\left(\frac{p^2-1+\sin^2\varphi}{p\sin^2\varphi}\right)^{-2b}.
\end{align}
Using \eqref{E:p2}, straightforward expansion of the right hand side of \eqref{E:ende} shows it to be equal to 
\[
\frac{b(2b-1)}{8}(1-L)^4+\frac{b(2b-1)}{4}(1-L)^5+(1-L)^4O(u^2(1-L)^2).
\]
\end{proof}

\begin{theorem}\label{T:asymp0} Let $q=e^a\in(0,1)$. Then
\begin{equation}\label{E:asymp0}
F(a,b,x)=P^b_{\mathbb U,e^{ix},1}(\gamma\subset A_q)\asymp\exp\left(\frac{b\pi x}{a}\right)
\end{equation}
as $a\nearrow0$, if either
\begin{equation}\label{E:condition1}
(b,x)\in([5/8,1]\cup[5/4,\infty))\times(0,\pi],
\end{equation}
or 
\begin{equation}\label{E:condition2}
(b,x)\in(1,5/4)\times(0,\pi)\text{ and }bx\le\pi.
\end{equation}
\end{theorem}

\begin{proof}
Equation \eqref{E:asymp0} holds as long as the difference between the upper and lower bound, estimated in Lemma \ref{L:upper-lower}, is not bigger than the lower bound from Corollary \ref{C:gcd}, i.e for $(b,x)\in[5/8,1]\times(0,\pi]$ or for $(b,x)\in[5/8,\infty)\times[0,\pi]$ such that $b x\le\pi$. 

For the remaining cases we use the following property of restriction measures: If $\gamma$ and $\gamma'$ are independent and with respective laws $P^b_{D,z_1,z_2}$ and $P^{b'}_{D,z_1,z_2}$, then the filling of $\gamma\cup\gamma'$ has law $P^{b+b'}_{D,z_1,z_2}$, see \cite{LSW:rest}. Here the filling of $\gamma\cup\gamma'$ is the smallest simply connected subdomain of $\overline{D}$ containing $\gamma\cup\gamma'$. This property and the definition extend to any finite number of hulls $\gamma_1,\dots,\gamma_n$ by induction.

Let now $b\in[5/4,\infty)$. Then $b=b_1+\cdots+b_n$ for some $n\in\mathbb Z^+$ and $b_1,\dots,b_n\in[5/8,1]$. Corollary \ref{C:gcd} provides a lower bound for $F(a,b,x)$ which is $\asymp \exp(b\pi x/a)$ as $a\nearrow0$ for all $x\in(0,\pi]$. For the upper bound, let $\gamma_1,\dots,\gamma_n$ be independent with respective laws $P^{b_k}_{\mathbb U,e^{ix},1}$, $1\le k\le n$. Then the event that the filling of $\gamma_1\cup\cdots\cup\gamma_n$ is contained in $A_q$ is a subset of the event $\{\gamma_1\cup\cdots\cup\gamma_n\subset A_q\}$. Hence
\[P^b_{\mathbb U, e^{ix},1}(\gamma\subset A_q)\le\prod_{k=1}^n P^{b_k}_{\mathbb U,e^{ix},1}(\gamma\subset A_q)\asymp e^{b\pi x/a}.\]
\end{proof}

\begin{remark} We believe the estimate \eqref{E:asymp0} holds for all $b\ge5/8$ and $x\in(0,\pi]$ but  a proof of this statement likely requires an upper bound closer to the lower bound than the upper bound we use.
\end{remark}


\begin{thebibliography}{99}

\bibitem{SLE83}
Bauer, Robert O., {\em Restricting SLE(8/3) to an annulus}, preprint, math.PR/0602391.

\bibitem{beffara}
Beffara, Vincent, Thesis, University of Paris, Orsay, 2003.

\bibitem{cardy:loop}
Cardy, John, {\it The O(n) model on the annulus}, preprint, math-ph/0604043.

\bibitem{hurwitz:1964}
Hurwitz, Adolf,
      {Vorlesungen \"uber allgemeine {F}unktionentheorie und
               elliptische {F}unktionen}, 4th edition, 
{Springer-Verlag}, {Berlin}, {1964}.

 \bibitem{ivanov}
 Ivanov, V. I. and Trubetskov, M. K., {\em Handbook of conformal mapping with computer-aided visualization}, CRC Press, Boca Raton, Florida, 1995.
 
\bibitem{lawler:2005}
Lawler, Gregory F., {Conformally invariant processes in the plane},
   Mathematical Surveys and Monographs,
   {114},
{American Mathematical Society},
     {Providence, RI},
       {2005}.

\bibitem{LSW:rest}
{Lawler, Gregory and Schramm, Oded and Werner, Wendelin},
      {\it Conformal restriction: the chordal case},
    {J. Amer. Math. Soc.},
  {\bf 16},
       {2003},
    {no. 4},
      {917--955 (electronic)}.

\bibitem{nehari:1952}
Nehari, Zeev, {Conformal mapping}, McGraw-Hill, New York, 1952.

\bibitem{werner:2005}
{Werner, Wendelin}, {\it The conformally invariant measure on self-avoiding loops}, preprint, math.PR/0511605.


\end{thebibliography}
\end{document}